\documentclass[12pt]{article}
\usepackage[centertags]{amsmath}
\usepackage{amsfonts}
\usepackage{amssymb}
\usepackage{amsthm}
\usepackage{newlfont}
\usepackage{graphicx}
\newtheorem{Theorem}{Theorem}[section]

\newtheorem{Lemma}{Lemma}[section]
\newtheorem{Remark}{Remark}[section]

\newtheorem{Example}{Example}
\theoremstyle{definition}

\theoremstyle{remark}

\numberwithin{equation}{section}

\begin{document}
\title{Dynamics of a rational multi-parameter second order difference equation with cubic numerator and quadratic monomial denominator}
\author{M. Shojaei\footnote{ Corresponding author, \ \ \textit{E-mail
address}: m\_shojaeiarani@aut.ac.ir,}}
\date{}

\maketitle
\begin{center}
\textit{\scriptsize Department of Applied
Mathematics, Faculty of Mathematics and Computer Science,
Amirkabir University of Technology, No. 424, Hafez Ave., Tehran,
Iran.}\\[0pt]
\end{center}
\begin{abstract}
The asymptotic behavior (such as convergence to an equilibrium,
convergence to a 2-cycle, and divergence to $\infty $) of solutions
of the following multi-parameter, rational, second order difference
equation $$x_{n+1}=\frac{ax_{n}^3+bx_{n}^2x_{n-1}+cx_{n}x_{n-1}^2+dx_{n-1}^3}{x_{n}^2}, \ \ \ x_{-1},x_{0}\in \Bbb{R},$$
is studied in this paper.\\
{\em Keywords}: Difference equation; equilibrium; 2-cycle; convergence; divergence
\end{abstract}
\section{Introduction}
Most of the work about rational difference equations treat the case where both numerator and denominator are linear polynomials. For second order rational difference equations with linear numerator and denominator we refer the reader to the monograph of Kulenovic and Ladass (\cite{KL}). In 2008, Sedaghat et al (\cite{DKMOS}) extended the existing results about second order rational difference equations to second order rational difference equations with quadratic numerator and linear denominator.

In this paper we extend the existing results to the following difference equation
\begin{equation}\label{formula2}
x_{n+1}=\frac{ax_{n}^3+bx_{n}^2x_{n-1}+cx_{n}x_{n-1}^2+dx_{n-1}^3}{x_{n}^2},
\end{equation}
which is a second order rational difference equation with cubic numerator and quadratic monomial denominator. The parameters $a,b,d$ are positive while the parameter $c$ and initial conditions $x_{-1},x_{0}$ could accept some negative values.

In (\cite{SH}) we investigated the dynamics of the
following difference equation
\begin{equation}\label{formula1}
x_{n+1}=\frac{ax_{n}^3+bx_{n}^2+cx_{n}+d}{x_{n}^3},
\end{equation}
where it was shown that in most cases every positive solution of
Eq.(\ref{formula1}) converges to either an equilibrium or, a
2-cycle.

In this part we study the asymptotic behavior of solutions of
Eq.(\ref{formula2}) including convergence to an equilibrium,
convergence to a 2-cycle, and divergence. Our analysis on the
dynamics of Eq.(\ref{formula2}) is essentially based on the dynamics
of Eq.(\ref{formula1}). The concepts of equilibrium point,
2-cycle, stability, asymptotic stability have been defined in the
first part and will not be repeated here. Moreover,
throughout the present paper we refer to some of the results in the
first part.

Divide both sides of Eq.(\ref{formula2}) by $x_{n}$ to obtain
$$\frac{x_{n+1}}{x_{n}}=a+b\left(\frac{x_{n-1}}{x_{n}}\right)+c\left(\frac{x_{n-1}}{x_{n}}\right)^{2}+d\left(\frac{x_{n-1}}{x_{n}}\right)^{3},$$
In the preceding equation substitute
$$t_{n}=\frac{x_{n}}{x_{n-1}},$$
to obtain
$$t_{n+1}=\frac{at_{n}^{3}+bt_{n}^{2}+ct_{n}+d}{t_{n}^{3}},$$
which simply is the first order Eq.(\ref{formula1}) (similar to the
first part we use the function $\phi (t)=(at^3+bt^2+ct+d)/t^3$,
which defines the right hand side of Eq.(\ref{formula1}), in the
present paper frequently). In fact the solutions of
Eq.(\ref{formula1}) are the successive ratios of the solutions of
Eq.(\ref{formula2}). So we call $\{t_{n}\}$ the sequence of ratios.
Eq.(\ref{formula2}) is a special semiconjugate factorization of
Eq.(\ref{formula1}) which is called semiconjugacy \ by \ ratios. For more about semiconjugacy and semiconjugacy by ratios
see \cite{Sedaghat2} and \cite{SHS} respectively. We analyze the
dynamics of Eq.(\ref{formula2}) using the dynamics of
Eq.(\ref{formula1}) which was studied in the first part.

Now we discuss about the initial conditions of Eq.(\ref{formula2}).
Since in this paper we studied the dynamics of positive solutions of
Eq.(\ref{formula1}) then the initial conditions of
Eq.(\ref{formula2}) should be chosen in such a way that the sequence
of ratios becomes positive eventually. If both $x_{-1}$ and $x_{0}$
are positive or negative then the sequence of ratios is positive
from the first step. On the other hand, if one of them be positive
and the other one be negative then the ratio $x_{n}/x_{n-1}$ may
never becomes positive or even the iteration process may stop. For
example if $\phi $ has a negative equilibrium which is attractive
then it attracts some ratios in a neighborhood around itself and
therefore such a ratio remains negative forever. Also, if at any
step the ratio equals zero then the iteration process stops. We
should avoid such cases. Although, the determination of these cases
in general is not possible but we are able to determine some of
them. Now, we mention one of them. Consider the function $\phi $ on
the interval $(-\infty ,0)$. Assume that $c_{-}<c\leq
-c^{*}=\sqrt{3bd}$ or, $c>-c^*$ and $\phi (x_{m})>0$ (note that
$x_{m}<0$ when $c>-c^*$). Then, there exists a unique number $r<0$
such that $\phi (r)=0$. Suppose that $r<r'<0$ is the unique number
such that $\phi (r')=r$. Then, it is evident that any ratio in
$(-\infty ,r)\cup (r',0)$ will eventually become positive after at
most three steps.
Note that if $\{x_{n}\}$ is a solution for Eq.(\ref{formula2}) then
so is $\{-x_{n}\}$. Also, first and third quadrants of
$\Bbb{R}^{2}$, namely $(0,\infty )^{2}$ and $(-\infty ,0)^{2}$, are
invariant under Eq.(\ref{formula2}). By the discussions in the
previous paragraph the ratio $x_{n}/x_{n-1}$ should be positive
eventually. Then there are two possibilities. Either $x_{n}>0$ or,
$x_{n}<0$ for all $n\geq n_{0}$ for some $n_{0}\in \Bbb{N}$. If the
second case occurs then the change of variable $y_{n}=-x_{n}$ (or
considering $\{-x_{n}\}$ as solution) reduces Eq.(\ref{formula2}) to
the first case. Therefore, without loss of generality we assume that
both of initial conditions are positive, hereafter.

In the first part we discussed (in great detail) about the
convergence of solutions of Eq.(\ref{formula1}) (or the sequence of
ratios) to both an equilibrium and a 2-cycle. Now, we want to study
the dynamics of solutions of Eq.(\ref{formula2}) in both of these
cases.

\section{Asymptotic stability when the sequence of ratios converges to an equilibrium}
\begin{Theorem}\label{Theorem8} Assume that the sequence
$\{x_{n}\}_{n=-1}^{\infty }$ is a positive solution for
Eq.(\ref{formula2}). Assume also that $\overline{t}$ is an
equilibrium of Eq.(\ref{formula1}) such that the sequence of ratios
$\{x_{n}/x_{n-1}\}_{n=0}^{\infty }$ converges to it.
\begin{description}
    \item[\it{(a)}] If $\overline{t}>1$ then $\{x_{n}\}$ diverges
    to $\infty $.
    \item[\it{(b)}] If $\overline{t}<1$ then $\{x_{n}\}$ converges to
    zero.
    \item[\it{(c)}] Assume that $\overline{t}=1$(or equivalently $a+b+c+d=1$). Let $\mathcal{S}=\{\phi ^{-n}(1)\}_{n=0}^{\infty
    }$. If $x_{0}/x_{-1}\in \mathcal{S}$ then $\{x_{n}\}$ is
    convergent to an equilibrium. Otherwise, $|b+2c+3d|\leq 1$ and also
    \begin{description}
        \item[\it{($c_{1}$)}] If $|b+2c+3d|<1$ then $\{x_{n}\}$
        converges to an equilibrium. Moreover, if $0<b+2c+3d<1$ then
        one of subsequences $\{x_{2n}\}$ and $\{x_{2n+1}\}$ will be
        increasing and the other one will be decreasing eventually
        while $\{x_{n}\}$ will be increasing or decreasing
        eventually if $-1<b+2c+3d\leq 0$.
        \item[\it{($c_{2}$)}] If $b+2c+3d=-1$ then $\{x_{n}\}$ will
        be increasing or decreasing eventually. In the later case
        $\{x_{n}\}$ converges to an equilibrium obviously. In the
        former case if $c>-3d$ then $\{x_{n}\}$ diverges to $\infty
        $.
        \item[\it{($c_{3}$)}] If $b+2c+3d=1$ then both of subsequences of even and odd
        terms will be increasing eventually. In particular, if $c>\frac{-2d}{a+d}-b$ then
        $\{x_{n}\}$ diverges to $\infty $.

    \end{description}
\end{description}
\end{Theorem}
Proof. (a) Since $\lim _{n\rightarrow \infty
}x_{n}/x_{n-1}=\overline{t}>1$ then there exist $L>1$ and $N\in
\Bbb{N}$ such that $x_{n}/x_{n-1}>L$ for all $n>N$. This simply
shows that $x_{n}\rightarrow \infty $ as $n\rightarrow \infty $. The
proof of (b) is similar and will be omitted.

(c) The equality $\overline{t}=1$ is simply equivalent to the
equality $a+b+c+d=1$. If $x_{0}/x_{-1}\in\mathcal{S}$ then there
exists $N\in \Bbb{N}$ such that $x_{n+1}/x_{n}=1$ for all $n>N$.
Thus $x_{n}$ remains constant for all $n> N$. Hence, $x_{n}$
converges to an equilibrium.

Next, assume that $x_{0}/x_{-1}\not \in \mathcal{S}$. Then, since
$\{x_{n}/x_{n-1}\}$ converges to the equilibrium $\overline{t}=1$ we
have $|\phi ^{'}(1)|\leq 1$ or equivalently
\begin{equation}\label{f11}
|b+2c+3d|\leq 1.
\end{equation}

On the other hand, since $a+b+c+d=1$ then we obtain by some
computations that
\begin{equation}\label{f12}
x_{n+1}-x_{n}=r_{n}(x_{n}-x_{n-1}), \ \ \  r_{n}=-
\left(b+c+d+\frac{c+d}{t_{n}}+\frac{d}{t_{n}^{2}}\right),
\end{equation}
notice that $r_{n}\rightarrow b+2c+3d$ as $n\rightarrow\infty $
since $t_{n}$ converges to 1. Therefore, by (\ref{f11}) there are
three cases to consider as follow:\\
Case \ I; $|b+2c+3d|<1$: Thus there exist $0<L<1$ and $N\in
\Bbb{N}$ such that $|r_{n}|<L$ for all $n> N$. So (\ref{f12})
implies for $n> N$ that
\begin{equation}\label{f13}
|x_{n+1}-x_{n}|<L|x_{n}-x_{n-1}|,
\end{equation}
thus we have (by induction) for $n\geq N$ that
$$|x_{n+1}-x_{n}|<L^{n-N}|x_{N+1}-x_{N}|,$$
Therefore
\begin{equation}\label{f14} \lim _{n\rightarrow \infty
}x_{n+1}-x_{n}=0.
\end{equation}

On the other hand we obtain from (\ref{f13}) for $n>N$ that
\begin{eqnarray*}
|x_{n+1}-x_{N}| & \leq & |x_{n+1}-x_{n}|+|x_{n}-x_{n-1}|+\ldots
  +|x_{N+1}-x_{N}|\\
   & < & \left(\sum _{i=0}^{n-N}L^{i}\right)|x_{N+1}-x_{N}|< \left(\sum _{i=0}^{\infty
   }L^{i}\right)|x_{N+1}-x_{N}|\\
  & = & \frac{|x_{N+1}-x_{N}|}{1-L},
\end{eqnarray*}
therefore, $\{x_{n}\}$ is bounded. This fact together with
(\ref{f14}) imply that $\{x_{n}\}$ is convergent.

On the other hand , (\ref{f12}) implies that
$t_{n}(t_{n+1}-1)=r_{n}(t_{n}-1)$. Therefore
\begin{eqnarray*}
  t_{n+1}t_{n}-1 &=& t_{n+1}t_{n}\mp t_{n}-1 \\
   &=& t_{n}(t_{n+1}-1)+(t_{n}-1)\\
   &=&(t_{n}-1)(r_{n}+1),
\end{eqnarray*}
or equivalently
\begin{equation}\label{fff13}
\frac{x_{n+1}}{x_{n-1}}-1=(t_{n}-1)(r_{n}+1),
\end{equation}

Note that since $|r_{n}|\rightarrow |b+2c+3d|<1$ as $n\rightarrow
\infty $ there exists $n_{0}\in \Bbb{N}$ such that $r_{n}+1>0$ for
all $n>n_{0}$. Also, we know that when $\phi '(1)=-(b+2c+3d)\in
(-1,0)$, $t_{n}$ oscillates alternately around $1$ while when
$-(b+2c+3d)\in [0,1)$, $t_{n}$ remains on one side of $1$ forever.
These facts together with (\ref{fff13}) complete the proof of
$(c_{1})$.

Case \ II; $b+2c+3d=-1$: In this case $\phi '(1)=\phi (1)=1$
(recall that this case occurs when $c=c_{m}$ and $1$ is the greater
equilibrium of $\phi$ or, $c=c_{M}$ and $1$ is the lower equilibrium
of $\phi $. This case also may occur when $x_{m}<1<x_{M}$)
. Therefore, it's evident that there exists an $n_{0}\in \Bbb{N}$
such that either $x_{n}/x_{n-1}<1$ or, $x_{n}/x_{n-1}>1$ for
$n>n_{0}$. If the former case occurs then $\{x_{n}\}$ is decreasing
for all $n>n_{0}$ and therefore it will converge to an equilibrium.

On the other hand if the later case occurs then $\{x_{n}\}$ is
increasing for $n>n_{0}$. Define the following function
\begin{equation*}
    r(t)=-\left( b+c+d+\frac{c+d}{t}+\frac{d}{t^2} \right), \ \ \ t>0,
\end{equation*}
note that $r(1)=1, \ r'(1)=c+3d$. Thus, if $c>-3d$ then $r'(1)>0$.
As a result, there exists $\epsilon >0$ such that $r(t)>1$ for all
$t\in (1,1+\epsilon )$. Therefore, since $r_{n}=r(t_{n})$,
$t_{n}=x_{n}/x_{n-1}>1$ for all $n>n_{0}$, and $t_{n}\rightarrow 1$
as $n\rightarrow \infty $ we conclude that there exists
$n_{1}>n_{0}$ such that $r_{n}>1$ for all $n>n_{1}$. Thus by
(\ref{f12}) the sequence of differences $\{x_{n+1}-x_{n}\}$ is
increasing for $n>n_{1}$.

This fact together with the fact that $\{x_{n}\}$ is increasing
eventually imply that $\{x_{n}\}$ diverges to $\infty $.

Case  III; $b+2c+3d=1$: In this case $\phi '(1)=-\phi
(1)=-1$ (recall that this case may occur when $c\geq c^{*}$ or,
$c<c^{*}$ and $1<x_{m}$. Also note that by Lemma 2(c) in \cite{SH}
this case never occurs when $x_{M}<1$). Therefore, it's evident that
there exists $n_{0}\in \Bbb{N}$ such that the sequence of ratios
oscillate alternately around $1$ for all $n>n_{0}$. Some
computations show that
\begin{equation}\label{f15}
x_{n+1}-x_{n-1}=\rho _{n}(x_{n}-x_{n-1}), \ \ \ \rho
_{n}=c+2d-\frac{c+d}{t_{n}}-\frac{d}{t_{n}^{2}},
\end{equation}
since $\rho _{n}=(t_{n}-1)[(c+2d)t_{n}+d]/t_{n}^2$, $c+2d=a>0$, and
$t_{n}\neq 1$ for all $n\geq 0$ then
\begin{equation}\label{ff16}
\rho _{n}(t_{n}-1)>0,
\end{equation}
where the equality $c+2d=a$ is gained by the subtraction of
equalities $a+b+c+d=b+2c+3d=1$. Now consider the consecutive ratios
$x_{2n}/x_{2n-1}$ and $x_{2n+1}/x_{2n}$ for $n>n_{0}$. Since these
ratios oscillate around $1$ alternately then one of them is greater
than $1$ and the other one is less than $1$. Without loss of
generality assume that $x_{2n+1}/x_{2n}<1<x_{2n}/x_{2n-1}$. Thus by
(\ref{ff16}) $\rho _{2n+1}<0<\rho _{2n}$. Therefore, (\ref{f15})
implies that $x_{2n+1}>x_{2n-1}$ and $x_{2n+2}>x_{2n}$, i.e., both
of subsequences of even and odd terms are increasing eventually.

Next, define the function
$$R(t)=r(t)r(\phi (t)), \ \ \ t>0,$$
where $r$ is defined in the previous case. Some algebra shows that
$$R(1)=1,\ \ \ R'(1)=0,\ \ \ R''(1)=2(b+c)(a+d)+4d,$$
Therefore, if $R''(1)>0$, i.e., $c>-2d/(a+d)-b$ then $1$ is a local
minimum point for $R$. As a result, there exists $\epsilon >0$ such
that $R(t)>1$ for all $t\in (1-\epsilon ,1+\epsilon ),t\neq 1$.
Therefore, since $r_{n+1}r_{n}=R(t_{n})$ and $t_{n}\rightarrow 1$ as
$n\rightarrow \infty $ then there exists $n_{1}>n_{0}$ such that
$r_{n+1}r_{n}>1$ for all $n>n_{1}$. Thus, we obtain from (\ref{f12})
that for $n>n_{1}$
$|x_{n+1}-x_{n}|=r_{n}r_{n-1}|x_{n-1}-x_{n-2}|>|x_{n-1}-x_{n-2}|$ or
equivalently, $|d_{n+1}|>|d_{n-1}|$ where $d_{n}$ is the sequence of
differences. Therefore, both of sequences $\{|d_{2n+1}|\}$ and
$\{|d_{2n}|\}$ are increasing. Hence, either they are convergent to
a positive number or divergent to $\infty $.

Finally, we claim that both of subsequences of even and odd terms
(and therefore $\{x_{n}\}$) diverge to $\infty $. Suppose for the
sake of contradiction that one of them is convergent or both of them
are convergent. If one of them is convergent and the other one is
divergent then this simply is a contradiction since the ratio
$x_{n}/x_{n-1}$ converges to $1$. On the other hand, if both of them
are convergent then by the same reason both of them should be
convergent to a same number. So $\{x_{n}\}$ is convergent. Thus,
$|d_{n+1}|=|x_{n+1}-x_{n}|\rightarrow 0$ as $n\rightarrow \infty $
which
is a contradiction. Therefore, $\{x_{n}\}$ diverges to $\infty $. The proof is complete.

\section{Asymptotic stability when the sequence of ratios converges to a 2-cycle}
\begin{Lemma}\label{Lemma5}\begin{description}
    \item[\it{(a)}] Eq.(\ref{formula1}) has a unique 2-cycle $(p,q)$ with
    $pq=1$ if and only if
    \begin{equation}\label{f16}
    \frac{a-c}{d}=\frac{d-b+1}{a}>2.
    \end{equation}
    \item[\it{(b)}] Assume that $(p,q)$ is an attractive 2-cycle of
    Eq.(\ref{formula1}) with $pq=1,p<1<q$. Then
    $$q+p\phi '(p)<0<p+q\phi '(q).$$
    \item[\it{(c)}] Assume that (\ref{f16}) holds. Then Eq.(\ref{formula2})
    has infinite number of 2-cycles. More precisely, the following set
    is the family of 2-cycles of Eq.(\ref{formula2})
    $$\mathcal{A}=\{(p',q')|\ p'/q'=p\ \   or\ \ p'/q'=q\}.$$
\end{description}
\end{Lemma}
Proof. (a) Assume that Eq.(\ref{formula1}) has a 2-cycle $(p,q)$
with $pq=1$. So $q^2=aq^3+bq^2+cq+d$ and $p^2=ap^3+bp^2+cp+d$.
Multiply the first equation by $p$ and the second equation by $q$
and apply some algebra to obtain
$$p+q=\frac{d-b+1}{a},$$
in a similar fashion multiply those two equations by $p^2$ and $q^2$
to obtain
$$p+q=\frac{a-c}{d}.$$

Therefore $(a-c)/d=(d-b+1)/a$. Since $pq=1$ and $p+q=(a-c)/d$ then
both $p$ and $q$ satisfy the following quadratic polynomial
\begin{equation}\label{f17}
X^2-\frac{a-c}{d}X+1=0,
\end{equation}
Therefore such a 2-cycle is unique. On the other hand, Eq.(\ref{f17}) should have positive determinant.
So $(a-c)/d>2$ and therefore (18) holds.

Next, suppose that (18) holds. Then, it's easy to verify that the
polynomial $G$ in Lemma 1 in \cite{SH} is factored by
Eq.(\ref{f17}). As a result, Eq.(\ref{formula1}) has a 2-cycle
$(p,q)$ with $pq=1$.

(b) At first we show that both of quantities $q+p\phi '(p)$ and
$p+q\phi '(q)$ have different signs. Since $(p,q)$ is an attractive
2-cycle of Eq.(\ref{formula1}) then
\begin{equation}\label{f18}
\phi '(p)\phi '(q)=(\phi ^2)'(p)\leq 1.
\end{equation}

On the other hand, since $pq=1$ then (a) implies that (\ref{f16})
holds. This fact together with
the fact that $p+q=(a-c)/d$ imply that
\begin{eqnarray}\label{f19}
\nonumber  p^2\phi '(p)+q^2\phi '(q) & =
&-\left(\frac{bq^2+2cq+3d}{q^2}+\frac{bp^2+2cp+3d}{p^2}\right)
  \\
 \nonumber & = & -\left(2b+2c(p+q)+3d((p+q)^2-2)\right) \\
 \nonumber & = & -\left(2b+2c\left(\frac{a-c}{d}\right)+3d\left[\left(\frac{a-c}{d}\right)^2-2\right]\right) \\
 \nonumber & = & -\left(2b+2a\frac{a-c}{d}+\frac{(a-c)^2}{d}-6d\right)\\
 \nonumber & = & -\left(2b+2(d-b+1)+\frac{(a-c)^2}{d}-6d\right)\\
 \nonumber & = & -\left(2+\frac{(a-c)^2-4d^2}{d}\right)\\
  &<&-2.
\end{eqnarray}

Thus (\ref{f18}), (\ref{f19}), and the equality $pq=1$ yield
$$(q+p\phi '(p))(p+q\phi '(q))=1+\phi '(p)\phi '(q)+p^2\phi '(p)+q^2\phi '(q)<2-2=0,$$
therefore, both of quantities $q+p\phi '(p)$ and $p+q\phi '(q)$ have
different signs. On the other hand, the equality $pq=1$ together
with (\ref{f17}) imply that
\begin{equation}\label{f20}
q+p\phi '(p)=q(1-b-2cq-3dq^2)=q(1-b+3d+(c-3a)q),
\end{equation}
similarly
\begin{equation}\label{f21}
p+q\phi '(q)=p(1-b+3d+(c-3a)p).
\end{equation}

If $p+q\phi '(q)<0<q+p\phi '(p)$ then (\ref{f20}) and (\ref{f21})
imply that
$$0<(1-b+3d+(c-3a)q)-(1-b+3d+(c-3a)p)=(c-3a)(q-p),$$
So since $p<q$ we obtain that $c>3a$ which simply contradicts
(\ref{f16}). Hence, $q+p\phi '(p)<0<p+q\phi '(q)$.

(d) The proof of (d) is clear and will be omitted.

The proof is complete.\\

The following theorem (whose proof somehow uses the ideas in the
proof of Theorem \ref{Theorem8}) discusses about the dynamics of
solutions of Eq.(\ref{formula2}) when the sequence of ratios converges to a 2-cycle.

\begin{Theorem}\label{Theorem9}Assume that the sequence
$\{x_{n}\}_{n=-1}^{\infty }$ is a positive solution for
Eq.(\ref{formula2}) and $(p,q)$ is a 2-cycle of Eq.(\ref{formula1})
such that the sequence of ratios $\{x_{n}/x_{n-1}\}_{n=0}^{\infty }$
converges to it.
\begin{description}
    \item[\it{(a)}] If $pq>1$ then $\{x_{n}\}$ diverges to $\infty $.
    \item[\it{(b)}] If $pq<1$ then $\{x_{n}\}$ converges to zero.
    \item[\it{(c)}] Assume that $pq=1$(or equivalently (\ref{f16}) holds).
    Let $\mathcal{S}=\{\phi ^{-n}(p),\phi  ^{-n}(q)\}_{n=0}^{\infty
    }$. If $x_{0}/x_{-1}\in \mathcal{S}$ then $\{x_{n}\}$ converges
    to a 2-cycle. Otherwise, we have $|\phi '(p)\phi '(q)|\leq 1$
    and we consider three cases as follow
    \begin{description}
        \item[$(c_{1})$] $|\phi '(p)\phi '(q)|<1$; In this case
    $\{x_{n}\}$ converges to a 2-cycle. Moreover, if $-1<\phi '(p)\phi '(q)<0$ then the
    subsequences $\{x_{4n}\}$ and $\{x_{4n+3}\}$ will be
    increasing and the other two will be decreasing eventually or vice
    versa while both of subsequences $\{x_{2n}\}$ and $\{x_{2n+1}\}$ will be
    increasing or decreasing eventually if $0\leq \phi '(p)\phi '(q)< 1$.
        \item[$(c_{2})$] $\phi '(p)\phi '(q)=1$; In this case both of
        subsequences $\{x_{2n}\}$ and $\{x_{2n+1}\}$ will be increasing or decreasing
        eventually. In the later case $\{x_{n}\}$ converges to a
        2-cycle. In the former case $\{x_{n}\}$ diverges to $\infty
        $ if
        $$p+q\phi '(q)+(\phi '(q))^2\phi ''(p)/2+\phi '(p)\phi ''(q)/2>0$$
        \item[$(c_{3})$] $\phi '(p)\phi '(q)=-1$; Let $l=-(p^2+(\phi '(q))^2q^2)+(q+p\phi '(p))\phi ''(q)/2+(p+q\phi '(q))(\phi '(q))^2$
        $\phi ''(p)/2$. If $l<0$ then all of subsequences $\{x_{4n}\},
        \{x_{4n+1}\},\{x_{4n+2}\},$ and $\{x_{4n+3}\}$ are
        decreasing eventually. In this case $\{x_{n}\}$ converges to a 2-cycle.
        If $l>0$ then all of subsequences $\{x_{4n}\},
        \{x_{4n+1}\},\{x_{4n+2}\},$ and $\{x_{4n+3}\}$ are
        increasing eventually. In this case $\{x_{n}\}$ diverges to $\infty $
        if
        $$-2s''(q)-2(s'(q))^2-s'(q)(\phi ^2)''(q)>0$$
        where
        $$s(t)=\frac{t\phi (t)\gamma (\phi (t))\theta (t)[\phi ^2(t)\theta (\phi ^2(t)+p]}{t\theta
        (t)+p},$$\\
         $$\gamma
         (t)=-\left(\frac{b}{pt}+\frac{c(t+p)}{p^2t^2}+\frac{d(t^2+pt+p^2)}{p^3t^3}\right),$$\\
         $$ \theta (t)=-\left(\frac{b}
         {qt}+\frac{c(t+q)}{q^2t^2}+\frac{d(t^2+qt+q^2)}{q^3t^3}\right).$$
\end{description}
\end{description}
\end{Theorem}
Proof. Throughout the proof we assume, without loss of generality,
that
\begin{equation}\label{f22}
t_{2n}\rightarrow p,\ \ \ t_{2n+1}\rightarrow q,\ \ \  \text{as} \
n\rightarrow \infty ,
\end{equation}
therefore
$$\frac{x_{n+2}}{x_{n}}=t_{n+2}t_{n+1}\rightarrow pq \ \ \ \text{as} \ n\rightarrow \infty ,$$
Thus if $pq>1$ then there exist $N\in \Bbb{N}$ and $L>1$ such that
$x_{n+2}/x_{n}>L$ for all $n>N$. This simply proves (a). In a
similar fashion (b) is proved. Now we proceed to (c). Since $pq=1$
then we assume, without loss of generality, that $p<1<q$ hereafter.
If $x_{0}/x_{-1}\in \mathcal{S}$ then there exists an integer $N$
such that $x_{N+1}/x_{N}=p$ or $x_{N+1}/x_{N}=q$. Thus
$x_{n+1}/x_{n}=p$ and $x_{n+2}/x_{n+1}=q$ for all $n\geq N$ or vice
versa. Therefore for $n\geq N$
$$\frac{x_{n+2}}{x_{n}}=pq=1,$$
which means that $\{x_{n}\}$ converges to a 2-cycle. Now assume that
$x_{0}/x_{-1}\not \in \mathcal{S}$. Then since the 2-cycle $(p,q)$
attracts the sequence of ratios $\{x_{n}/x_{n-1}\}$ we have $|\phi
'(p)\phi '(q)|\leq 1$.

$(c_{1})$ Define $D_{n}=x_{n}-x_{n-2}$. Therefore, using (\ref{f22})
and hopital law in calculous one can write
\begin{eqnarray*}
  \lim _{n\rightarrow \infty }\left|\frac{D_{2n+2}}{D_{2n}}\right| &=& \lim _{n\rightarrow \infty }\left|\frac{t_{2n+2}t_{2n+1}t_{2n}t_{2n-1}-t_{2n}t_{2n-1}}{t_{2n}t_{2n-1}-1}\right|  \\
   &=&\lim _{t\rightarrow q} \left|\frac{t\phi (t)\phi ^2(t)\phi ^3(t)-t\phi (t)}{t\phi
   (t)-1}\right|\\
   & = & |\phi '(p)\phi '(q)|\\
   & < & 1.
\end{eqnarray*}

In a similar fashion
$$\lim _{n\rightarrow \infty }\left|\frac{D_{2n+1}}{D_{2n-1}}\right|=|\phi '(p)\phi '(q)|<1, $$
Consequently, there exist $n_{0}\in \Bbb{N}$ and $0<L<1$ such that
for $n>n_{0}$
$$|D_{2n+2}|<L|D_{2n}|,\ \ \ |D_{2n+1}|<L|D_{2n-1}|.$$

Therefore, by an analysis precisely similar to what was applied in
Theorem \ref{Theorem8}(c) it could be shown that both of
subsequences $\{x_{2n}\}$ and $\{x_{2n+1}\}$ are convergent and
hence $\{x_{n}\}$ converges to a 2-cycle.

Some calculations show that
\begin{equation}\label{f23}
t_{n+1}-q=\gamma _{n}(t_{n}-p), \ \ \ \gamma
_{n}=-\left(\frac{b}{pt_{n}}+\frac{c(t_{n}+p)}{p^2t_{n}^2}+\frac{d(t_{n}^2+pt_{n}+p^2)}{p^3t_{n}^3}\right),
\end{equation}
and
\begin{equation}\label{f24}
t_{n+1}-p=\theta _{n}(t_{n}-q), \ \ \ \theta
_{n}=-\left(\frac{b}{qt_{n}}+\frac{c(t_{n}+q)}{q^2t_{n}^2}+\frac{d(t_{n}^2+qt_{n}+q^2)}{q^3t_{n}^3}\right).
\end{equation}

by (\ref{f22}) we obtain
\begin{equation}\label{f25}
\gamma _{2n} \rightarrow \phi '(p), \ \ \ \theta _{2n+1}\rightarrow
\phi '(q),\ \ \ \text{as} \ n\rightarrow \infty .
\end{equation}

Now, suppose that $-1<\phi '(p)\phi '(q)< 0$. By Lemma
\ref{Lemma5}(b), $\phi '(p)<-q/p<0$. So $\phi '(q)>0$. Therefore, in a neighborhood around $p$ and another neighborhood around $q$ $\phi$ is decreasing and increasing respectively. This fact together with (\ref{f22}) imply that there exists $n_{0}\in BBb{N}$ such that for
$n\geq n_{0}$
\begin{description}
    \item[\it{(i)}] either $t_{4n}<p,t_{4n+1}>q,t_{4n+2}>p,t_{4n+3}<q$ or,
    \item[\it{(ii)}] $t_{4n}>p,t_{4n+1}<q,t_{4n+2}<p,t_{4n+3}>q.$
\end{description}
On the other hand, (\ref{f23}) and (\ref{f24}) imply that
\begin{eqnarray*}
  t_{n+4}t_{n+3}t_{n+2}t_{n+1}-1 &=&  t_{n+4}t_{n+3}t_{n+2}t_{n+1}\mp p^2t_{n+3}t_{n+1}-p^2q^2 \\
   &=& t_{n+3}t_{n+1}(t_{n+4}t_{n+2}-p^2)+p^2(t_{n+3}t_{n+1}-q^2) \\
   &=& t_{n+3}t_{n+1}(t_{n+4}t_{n+2}\mp pt_{n+2}-p^2)+p^2(t_{n+3}t_{n+1}\mp qt_{n+1}-q^2) \\
   &=& t_{n+3}t_{n+1}[t_{n+2}(t_{n+4}-p)+p(t_{n+2}-p)]+ \\
   & & p^2[t_{n+1}(t_{n+3}-q)+q(t_{n+1}-q)]\\
   &=& t_{n+3}t_{n+1}(t_{n+2}\theta _{n+3}\gamma  _{n+2}+p)(t_{n+2}-p)+ \\
   & &   p^2(t_{n+1}\gamma _{n+2}\theta  _{n+1}+q)(t_{n+1}-q)\\
   &=& [t_{n+3}t_{n+1}\theta _{n+1}(t_{n+2}\theta _{n+3}\gamma  _{n+2}+p)+p^2(t_{n+1}\gamma _{n+2}\theta  _{n+1}+q)]\times \\
   & & (t_{n+1}-q),
\end{eqnarray*}
Therefore
\begin{equation}\label{f27}
\frac{x_{n+4}}{x_{n}}-1=\lambda _{n}(t_{n+1}-q), \  \ \lambda
_{n}=\& t_{n+3}t_{n+1}\theta _{n+1}(t_{n+2}\theta _{n+3}\gamma
 _{n+2}+p)+p^2(t_{n+1}\gamma _{n+2}\theta
 _{n+1}+q),
\end{equation}
In a similar fashion one can write
\begin{equation}\label{f28}
\frac{x_{n+4}}{x_{n}}-1=\xi _{n}(t_{n+1}-p), \ \ \ \xi
_{n}=\&t_{n+3}t_{n+1}\gamma _{n+1}(t_{n+2}\gamma _{n+3}\theta
 _{n+2}+q)+q^2(t_{n+1}\theta _{n+2}\gamma_{n+1}+p),
\end{equation}
notice that (\ref{f22}) and (\ref{f25}) imply that
\begin{equation}\label{f29}
\lambda _{2n}\rightarrow (\phi '(p)\phi '(q)+1)(p+q\phi '(q)), \ \ \
\xi _{2n+1}\rightarrow (\phi '(p)\phi '(q)+1)(q+p\phi '(p)), \ \ \
\text{as}\ n\rightarrow \infty .
\end{equation}

Consequently, by the fact that $\phi '(p)\phi '(q)>-1$, Lemma
\ref{Lemma5}(b), (\ref{f27}), (\ref{f28}), and (\ref{f29}) the
subsequences $\{x_{4n}\}$ and $\{x_{4n+3}\}$ will be increasing
while the other two will be decreasing eventually if $(i)$
holds. Otherwise, the subsequences $\{x_{4n}\}$ and $\{x_{4n+3}\}$
will be decreasing while the other two will be increasing
eventually.

Next, assume that $0\leq \phi '(p)\phi '(q)<1$. By Lemma
\ref{Lemma5}(b) $\phi '(p)<-q/p<0$. Thus $\phi '(q)\leq 0$. If $\phi
'(q)<0$ then in a neighborhood around $p$ and another neighborhood
around $q$ $\phi $ is decreasing. As a result by (\ref{f22}) we
conclude that there exists an integer $n_{0}\in \Bbb{N}$ such that
for $n\geq n_{0}$
\begin{description}
    \item[(i)] either  $t_{2n}>p,t_{2n+1}<q$ or,
    \item[(ii)] $t_{2n}<p,t_{2n+1}>q$.
\end{description}

If, on the other hand $\phi '(q)=0$ then $q=x_{m}$ or $q=x_{M}$.
It's easy to show that if $q=x_{m}$ then case $(i)$ occurs
while case $(ii)$ occurs if $q=x_{M}$. By an analysis somehow
similar to that of applied for the expression $x_{n+4}/x_{n}-1$ we
obtain
\begin{equation}\label{f31}
\frac{x_{n+2}}{x_{n}}-1=\lambda '_{n}(t_{n+1}-q)=\xi
'_{n}(t_{n+1}-p), \ \ \ \lambda '_{n}=t_{n+1}\theta _{n+1}+p, \ \ \
\xi '_{n}=t_{n+1}\gamma _{n+1}+q,
\end{equation}
with
\begin{equation}\label{f32}
\lambda '_{2n}\rightarrow p+q\phi '(q),\ \ \ \xi '_{2n+1}\rightarrow
q+p\phi '(p), \ \ \ \text{as} \ n\rightarrow \infty .
\end{equation}

Therefore, Lemma \ref{Lemma5}(b), (\ref{f31}), and (\ref{f32}) imply
that both of subsequences $\{x_{2n}\}$ and $\{x_{2n+1}\}$ are
decreasing eventually if $(i)$ holds and vice versa if
$(ii)$ holds.

$(c_{2})$ By an analysis precisely similar to what was applied for
the case $0\leq \phi '(p)\phi '(q)<1$ in $(c_{1})$ one can prove
that both of subsequences $\{x_{2n}\}$ and $\{x_{2n+1}\}$ are
increasing or decreasing eventually. If the later case occurs (note
that this case occurs when case $(i)$ in $(c_{1})$ occurs,
i.e., $t_{2n}>p,t_{2n+1}<q$ for $n>n_{0}$) then $\{x_{n}\}$
converges to a 2-cycle obviously.
Now assume that the former case occurs (note that in this case
$t_{2n}<p,t_{2n+1}>q$ for $n>n_{0}$). Then
\begin{eqnarray*}
\nonumber  \frac{D_{n+2}}{D_{n}} &=& \frac{x_{n+2}-x_{n}}{x_{n}-x_{n-2}}
\\
\nonumber   &=& \frac{t_{n}t_{n-1}(t_{n+2}t_{n+1}-1)}{t_{n}t_{n-1}-1}, \\
\end{eqnarray*}
Therefore (\ref{f31}), (\ref{f23}), (\ref{f24}), and some algebra
imply that
\begin{equation}\label{f33}
D_{n+2}=s_{n}D_{n}, \ \ \ s_{n}=\frac{t_{n}t_{n-1}\gamma _{n}\theta
_{n-1}\lambda '_{n}}{\lambda '_{n-2}}.
\end{equation}

Some computations show that
\begin{equation}\label{f34}
\gamma (p)=\phi '(p), \ \ \ \theta (q)=\phi '(q), \ \ \ \gamma
'(p)=\frac{\phi ''(p)}{2}, \ \ \ \theta '(q)=\frac{\phi ''(q)}{2},
\end{equation}
Notice that $s_{n}=s(t_{n-1})$ and by (\ref{f34}) $s(q)=1$. Also
using (\ref{f34}) and some algebra we obtain that
$$s'(q)=p+q\phi '(q)+(\phi '(q))^2\frac{\phi ''(p)}{2}+\phi '(p)\frac{\phi ''(q)}{2}>0,$$
As a result there exists $\epsilon >0$ such that $s(t)>1$ for all
$t\in (q,q+\epsilon )$. Therefore, since $s_{2n}=s(t_{2n-1})$,
$t_{2n-1}>q$ for all $n>n_{0}$, and $t_{2n-1}\rightarrow q$ as
$n\rightarrow \infty $ then there exists $n_{1}>n_{0}$ such that
$s_{2n}>1$ for all $n>n_{1}$. Thus by (\ref{f33}) the sequence
$\{D_{2n}\}$ is increasing eventually.

Consequently, since $\{x_{2n}\}$ is increasing eventually then
$\{x_{2n}\}$ should be divergent to $\infty $ and hence by
(\ref{f22}) $\{x_{2n+1}\}$ should be divergent to $\infty$, too.
This means that $\{x_{n}\}$ is divergent to $\infty $.

$(c_{3})$ Note that similar to the case $-1\leq \phi '(p)\phi
'(q)<0$ in $(c_{1})$ there exists $n_{0}\in \Bbb{N}$ such that
either $t_{4n}<p,t_{4n+1}>q,t_{4n+2}>p,t_{4n+3}<q$ or,
$t_{4n}>p,t_{4n+1}<q,t_{4n+2}<p,t_{4n+3}>q$ for $n>n_{0}$. Consider
the quantities $\lambda _{n}$ and $\xi _{n}$ in $(c_{1})$ and define
the following functions for $t>0$
\begin{eqnarray*}
  \lambda (t)=t\phi ^2(t)\theta (t)[\phi (t)\theta (\phi ^2(t))\gamma (\phi (t))+p]+p^2[t\gamma (\phi (t))\theta (t)+q],\\
  \xi (t)=t\phi ^2(t)\gamma (t)[\phi (t)\gamma (\phi ^2(t))\theta (\phi (t))+p]+p^2[t\theta (\phi (t))\gamma
  (t)+q],
\end{eqnarray*}
Notice that $\lambda _{n}=\lambda (t_{n+1}), \xi _{n}=\xi
(t_{n+1})$, and by (\ref{f34}) $\lambda (q)=\xi (p)=0$. Also by
(\ref{f34}) and some algebra we have
$$(\phi '(q))^2\xi '(p)=\lambda '(q)=l,$$
Therefore both of quantities $\xi '(p)$ and $\lambda '(q)$ have the
same signum. Now assume that both of them are negative, i.e., $l<0$.
Then there are neighborhoods around $p$ and $q$ that $\xi $ and
$\lambda $ are decreasing on respectively. Assume that
$t_{4n}<p,t_{4n+1}>q,t_{4n+2}>p,t_{4n+3}<q$ for $n>n_{0}$. Thus
since $\lambda _{4n}=\lambda (t_{4n+1}),\lambda _{4n+2}=\lambda
(t_{4n+3}), \xi _{4n+1}=\xi (t_{4n+2})$, and $\xi _{4n+3}=\xi
(t_{4n+4})$ then by (\ref{f22}) we obtain that there exists
$n_{1}>n_{0}$ such that for $n>n_{1}$
$$\lambda _{4n}<0<\lambda _{4n+2}, \ \ \ \xi _{4n+1}<0<\xi _{4n+3},$$
Consequently by (\ref{f27}) and (\ref{f28}) we conclude that all of
subsequences $\{x_{4n}\}$,$\{x_{4n+1}\}$,\\
$\{x_{4n+2}\},$ and $\{x_{4n+3}\}$ are decreasing eventually (note that similar result
obtains if $t_{4n}>p,t_{4n+1}<q,t_{4n+2}<p,t_{4n+3}>q$ for $n>n_{0}$). As a result, all of these four subsequences
are convergent and by the fact that $x_{n+2}/x_{n}\rightarrow 1$ as $n\rightarrow \infty $ we obtain that both of subsequences $\{x_{4n}\},\{x_{4n+2}\}$ of even terms should be convergent to a same number. The same result holds for the subsequences $\{x_{4n+1}\},\{x_{4n+3}\}$ of odd terms. Hence, $\{x_{n}\}$ converges to a 2-cycle.

Next, suppose that $l>0$. Then similar arguments show that all of
subsequences $\{x_{4n}\}$,$\{x_{4n+1}\}$,$\{x_{4n+2}\}$, and $\{x_{4n+3}\}$ are increasing eventually. Define
the function $$S(t)=s(t)s(\phi ^2(t)), \ \ \ t>0,$$ using the fact
that $s(q)=\phi '(p)\phi '(q)=-1$ and by some algebra we obtain that
$$S(q)=1, \ \ \ S'(q)=0,\ \ \ S''(q)=-2s''(q)-2(s'(q))^2-s'(q)(\phi ^2)''(q)>0,$$
Thus $q$ is a local minimum point for $S$. So there exists $\epsilon
>0$ such that $S(t)>1$ for $t\in (q-\epsilon ,q+\epsilon ),t\neq q$.
Therefore, since $s_{4n+2}s_{4n}=s(t_{4n+1})s(t_{4n-1})=S(t_{4n-1})$
and $t_{4n-1}\rightarrow q$ as $n\rightarrow \infty $ then there
exists $n_{0}\in \Bbb{N}$ such that $s_{4n+2}s_{4n}>1$ for all
$n>n_{0}$. As a result (\ref{f33}) implies that
$|D_{4n+4}|=s_{4n+2}s_{4n}|D_{4n}|>|D_{4n}|$, i.e., the sequence
$\{|D_{4n}|\}$ is increasing eventually. Thus, either it converges
to a positive number or, diverges to $\infty $.

We claim that both of subsequences of even terms, i.e., $\{x_{4n}\}$
and $\{x_{4n+2}\}$ are divergent to $\infty $ (and therefore by
(\ref{f22}) the other two subsequences are divergent, too. Hence,
$\{x_{n}\}$ diverges to $\infty $). Otherwise, at least one of them
should be convergent and therefore since $x_{4n+2}/x_{4n}\rightarrow
1$ as $n\rightarrow\infty $ we conclude that both of them are
convergent. As a result $D_{4n}\rightarrow 0$ as $n\rightarrow
\infty $ which
simply is a contradiction. The proof is complete.

\begin{Remark}\label{Remark5} In Theorem \ref{Theorem8} and Theorem
\ref{Theorem9} dynamical behavior of solutions of
Eq.(\ref{formula2}) was studied where the sequence of ratios
converges to an equilibrium and a 2-cycle respectively. By Theorem 4
and Theorem 5 in \cite{SH} we know that one of these two cases occur
definitely when $c\geq c^*$ or, $c<c^*,x_{M} \leq \overline{t}$ or,
$c<c^*,x_{m}\leq \overline{t}\leq x_{M}$. But if $c<c^*$ and
$\overline{t}<x_{m}$ the sequence of ratios may fail to be
convergent to an equilibrium or a 2-cycle. In this case according to
Theorem $6(a)$ in \cite{SH} the interval $I=[\phi (x_{m}),\phi
^2(x_{m})]$ is invariant under hypothesis \emph{(H)} or even ratios
eventually end up in $I$ if $c\leq c_{1}^*$. Therefore, if
$x_{0}/x_{-1}\in I$ and \emph{(H)} holds, or $x_{0}/x_{-1}\not \in
I$ but $c\leq c_{1}^*$ then $\{x_{n}\}$ diverges to $\infty $ when
$\phi (x_{m})\geq 1$ while $\{x_{n}\}$ converges to zero when $\phi
^2(x_{m})\leq 1$ obviously.
\end{Remark}

\section{Some examples}
\begin{Example}\label{Ex1} Consider the first example in Remark 4 in
\cite{SH}. Note that $c>c_{-}\approx -4.1305$ where $c_{-}$ is the
unique negative root of the cubic polynomial $Q$ in Theorem 1 in
\cite{SH}. So By Theorem $1(b)$ in \cite{SH} nonpositive iterations
of Eq.(\ref{formula1}) do not occur. In this example
Eq.(\ref{formula1}) has two equilibria and no 2-cycle. Some
computations show that
$$x_{m}\approx 0.7133,\overline{t}_{1}\approx 0.7845,\overline{t}_{2}=1,\delta \approx0.5833,$$
where $\delta $ has been defined in Theorem $7(a)$ in \cite{SH}.
Notice that here $a+b+c+d=1$ and hence one of two equilibria is $1$.
By Theorem $7(a_{1})$ if $t_{0}\in (\delta ,\overline{t}_{2})$ then
$\{t_{n}\}$ converges to $\overline{t}_{1}$ otherwise, it converges
to $\overline{t}_{2}$. Moreover, if $t_{0}\not\in (\delta
,\overline{t}_{2})$ then $\{t_{n}\}$ converges to $t_{2}$ from the
right.

Therefore, since $\overline{t}_{1}<1$, $a+b+c+d=1$, $b+2c+3d=-1$,
and $c>-3d$ then Theorem \ref{Theorem8}($c_{2}$) and Theorem
\ref{Theorem8}(b) imply that
\begin{description}
    \item[\it{(i)}] If $x_{0}/x_{-1}\in (\delta ,\overline{t}_{2})$ then
    $\{x_{n}\}$ converges to zero.
    \item[\it{(ii)}] If $t_{0}\not\in (\delta
,\overline{t}_{2})$ then $\{x_{n}\}$ diverges to $\infty
$.
\end{description}
\end{Example}

\begin{Example}\label{Ex2} In Eq.(\ref{formula2}) set
$a=0.2,b=1.7,c=-2,d=1.1$. So $c>c_{-}\approx -2.8540$ and therefore
similar to the arguments in the previous example nonpositive
iterations of Eq.(\ref{formula1}) do not occur. In this example
Eq.(\ref{formula1}) has a unique equilibrium $\overline{t}=1$
(notice that $a+b+c+d=1$) and two 2-cycles $(p_{1},q_{1})\approx
(0.2262,63.6517)$ and $(p_{2},q_{2})\approx (0.5110,4.1111)$. Since
$c>c^*=-\sqrt{3bd}\approx -2.3685$ then by Theorem $4(c)$ in
\cite{SH} $\{t_{n}\}$ converges to $1$ if $t_{0}\in (p_{2},q_{2})$
and converges to the 2-cycle $(p_{1},q_{1})$ if $t_{0}\in
(0,p_{2})\cup (q_{2},\infty )$.

Therefore, since $p_{1}q_{1}>1$, $a+b+c+d=b+2c+3d=1$, and
$c>\frac{-2d}{a+d}-b$ then Theorem \ref{Theorem8}($c_{3}$), Theorem
\ref{Theorem9}(a), and Theorem \ref{Theorem9}(c) imply that
\begin{description}
    \item[\it{(i)}] If $x_{0}/x_{-1}\in (0,\infty )\setminus
    \{p_{1},p_{2},\overline{t},q_{2},q_{1}\}$ then $\{x_{n}\}$
    diverges to $\infty $.
    \item[\it{(ii)}] if $x_{0}/x_{-1}=\overline{t}$ then
    $\{x_{n}\}$ converges to an equilibrium.
    \item[\it{(iii)}] If $x_{0}/x_{-1}\in
    \{p_{1},p_{2},q_{2},q_{1}\}$ then $\{x_{n}\}$ converges to a
    2-cycle.
\end{description}
\end{Example}

\begin{Example}\label{Ex3} In Eq.(\ref{formula2}) set
$a=0.1,b=1.79,c=-2,d=1$. Thus again $c>c_{-}\approx-2.7295$ which
similar to the previous examples this guarantees that all iterations
of Eq.(\ref{formula1}) remain positive forever. In this example
Eq.(\ref{formula1}) has a unique equilibrium $\overline{t}\approx
0.9423$ and three 2-cycles $(p_{1},q_{1})\approx (0.1024,759.2585)$,
$(p_{2},q_{2})=(0.6021,2.1370)$ and $(p_{3},q_{3})\approx
(0.7298,1.3702)$. Since $c>c^*\approx -5.37$ then by Theorem $4(d)$
in \cite{SH} $\{t_{n}\}$ converges to the 2-cycle $(p_{1},q_{1})$ if
$t_{0}\in (0,p_{2})\cup (q_{2},\infty )$ and converges to the
2-cycle $(p_{3},q_{3})$ if $t_{0}\in (p_{2},q_{2})\setminus
\{\overline{t}\}$. Notice that $p_{3}q_{3}=1$. This is evident since
in this example (\ref{f16}) holds easily.

Consequently, since $p_{3}q_{3}=1$, $0<\phi' (p_{3})\phi
'(q_{3})<1$, and $p_{1}q_{1}>1$ then by Theorem \ref{Theorem8},
Theorem \ref{Theorem9}(a), and Theorem \ref{Theorem9}$(c_{1})$ we
conclude that
\begin{description}
    \item[\it{(i)}] If $x_{0}/x_{-1}\in (0,p_{2})\cup (q_{2},\infty
    )$  then $\{x_{n}\}$ diverges to $\infty $.
    \item[\it{(ii)}] If $x_{0}/x_{-1}\in [p_{2},q_{2}]\setminus \{\overline{t}\}$ then
    $\{x_{n}\}$ converges to a 2-cycle.
    \item[\it{(iii)}] If $x_{0}/x_{-1}=\overline{t}$ then
    $\{x_{n}\}$ simply converges to an equilibrium.
\end{description}
\end{Example}

\end{document}